\documentclass[12pt]{article}

\usepackage{latexsym}
\usepackage{amsmath}
\usepackage{amsfonts}
\usepackage{amssymb}

\newcommand{\Res}{\operatorname{Res}}

\newcommand{\bean}{\begin{eqnarray}}
\newcommand{\eean}{\end{eqnarray}}
\newcommand{\bea}{\begin{eqnarray*}}
\newcommand{\eea}{\end{eqnarray*}}
\newcommand{\bsa}{\begin{subarray}{c}}
\newcommand{\esa}{\end{subarray}}
\newcommand{\bi}{\begin{itemize}}
\newcommand{\ei}{\end{itemize}}

\newtheorem{lemma}{Lemma}[section]
\newtheorem{thm}[lemma]{Theorem}

\newtheorem{propn}[lemma]{Proposition}

\begin{document}

\title{\bf Remarks on the McKay Conjecture}
\author{Geoffrey Mason\thanks{Supported by grants from the NSF, NSA,  and faculty
research funds granted by the University of
California at Santa Cruz.}  \\
    Department of Mathematics \\
     University of California at Santa Cruz \\
     Ca. 95064, USA \\
   gem@cats.ucsc.edu }
\date{}
\maketitle
 \begin{abstract}
\noindent
The McKay Conjecture (MC) asserts the existence of a bijection between the (inequivalent) complex irreducible representations of degree coprime to $p$ ($p$ a prime) of a finite group $G$ and those of the subgroup $N$, the normalizer of Sylow $p$-subgroup. In this paper we observe that MC implies the existence of analogous bijections involving various pairs of algebras,  including
certain crossed products, and that MC is \emph{equivalent} to the analogous
statement for (twisted) quantum doubles. Using standard conjectures in orbifold conformal field theory, MC is \emph{equivalent} to parallel statements about holomorphic orbifolds
$V^G, V^N$. There is a uniform formulation of MC covering these different situations which involves
quantum dimensions of objects in pairs of ribbon fusion categories. 

\end{abstract}

\section{Introduction}
The following notation will be fixed throughout the paper:
 $G$ is a finite group, $p$ a prime, $P$ a Sylow $p$-subgroup of $G$, $N = N(P)$ the \emph{normalizer} of $P$ in $G$, $X$ a (finite, non-empty, left-) $G$-set. All algebras and modules are finite-dimensional and defined over $\mathbb{C}$. $\mathbb{C}[G]$ is the  group algebra of $G$ and $\mathbb{C}[G]^*$ the \emph{dual} group algebra. 

\medskip
For an algebra $A$, let
\begin{eqnarray}\label{Mpairdef}
\mu(A) = \# \mbox{inequivalent simple $A$-modules of dimension coprime to $p$}.
\end{eqnarray}
We say that a pair of algebras $(A, B)$ is an \emph{M-pair} in case $\mu(A) = \mu(B)$. The McKay Conjecture (MC) is the assertion that
$(\mathbb{C}[G], \mathbb{C}[N])$ is an M-pair. The reader may consult the paper \cite{IMN} of Isaacs, Malle and Navarro for the current status of this conjecture. The idea of the present paper is to \emph{extend}  
 MC beyond its original formulation for groups. First we show how it may be extended to large classes of algebras that are not group algebras.  Examples include \emph{crossed product} algebras, where we show that $(\mathbb{C}[H]^* \#_{\sigma} \mathbb{C}[G], \mathbb{C}[H]^* \#_{\sigma} \mathbb{C}[N])$ is an M-pair. Here, $G$ acts on the group $H$ and $\sigma$ is a certain $2$-cocycle.  (See \cite{Mo}, \cite{KMM} for background.) A particularly interesting case is that of quantum doubles $D(G)$ (see \cite{D}, \cite{M1} and below for more details). In this case we establish
\begin{eqnarray*}\label{AMCequiv1}
\mbox{MC is true if, and only if, $(D(G), D(N))$ is an M-pair for all $G$ and $N$.}
\end{eqnarray*}
Note that quantum doubles $D(G)$ are generally not group algebras (unless $G$ is abelian).

\medskip
For a multiplicative $3$-cocycle $\omega \in Z^3(G, \mathbb{C}^*)$, we show that 
MC implies the same result for \emph{twisted quantum doubles}. That is,  $(D^{\omega}(G), D^{\omega}(N))$ is an M-pair. Now there is a standard Ansatz in orbifold conformal field theory (CFT) due to Dijkgraaf-Pasquier-Roche \cite{DPR} which, when interpreted appropriately, says that the tensor category $D^{\omega}(G)$-Mod is  equivalent to the module category $V^G$-Mod 
of a so-called  holomorphic $G$-orbifold for a suitable
vertex operator algebra $V$ admitting $G$ as automorphisms. Therefore, granted the DPR conjecture, MC is  \emph{equivalent} 
to a CFT-formulation involving a bijection between certain sets of simple modules for $V^G$ and $V^N$. It is not necessary for the reader to be familiar with this language; the point  is simply that modules for
$V^G$ are infinite-dimensional and the idea of an M-pair based on definition (\ref{Mpairdef}) makes no sense. In fact, all three types of M-pairs that we have discussed (i.e. for groups, (quasi-)Hopf algebras and orbifolds) may be uniformly described in the following setting:  a pair of ribbon categories admitting a bijection between objects
whose quantum dimension is integral and coprime to $p$. 

\medskip
All of the proofs in this paper are elementary and involve nothing beyond a few facts about finite groups, their representations, and their cohomology. In Section 2 we discuss some algebras 
$D_X(G)$ constructed from $G$ and a $G$-set $X$ and show that MC implies that $(D_X(G), D_X(N)$
is an M-pair. We also establish (\ref{AMCequiv1}). In Section 3 we carry out the twisted analog of this construction. Together, these results 
cover several of the connections with crossed products and twisted quantum doubles mentioned above.  In Section 4 we discuss the connections with CFT and ribbon categories. We assure the reader that no knowledge of CFT is required to understand the contents of this paper.

\medskip
We thank Siu-Hung Ng, Robert Boltje and Gabriel Navarro for their interest and comments on earlier versions of this paper.

\section{The algebras $D_X(G)$}
We use the following additional notation: for $H \leq G, g \in G, H^g = \{g^{-1}hg \ | \ h \in H \}$. For $x \in X$,  Stab$_G(x) = \{g \in G | \ g.x = x \}$.

\bigskip
We now introduce the algebras $D_X(G)$, which were mentioned briefly in \cite{M1}. 
Let $\mathbb{C}[X]^*$ be the space of complex-valued functions on $X$. One sees that it is a $G$-module algebra, as follows.
The algebra structure is pointwise multiplication, with basis the Dirac delta functions
\begin{eqnarray*}
e(x): y \mapsto \delta_{x, y} \ \ x, y \in X.
\end{eqnarray*}
Thus
\begin{eqnarray*}
e(x)e(y) = \delta_{x, y}e(x).
\end{eqnarray*}
$G$ acts on the left of $\mathbb{C}[X]^*$ as algebra automorphisms via
\begin{eqnarray*}
g: e(x) \mapsto e(g.x).
\end{eqnarray*}
Consider the linear space
\begin{eqnarray}\label{DXGdef}
D_X(G) = \mathbb{C}[X]^* \otimes_{\mathbb{C}} \mathbb{C}[G].
\end{eqnarray}
It becomes an algebra via the product 
\begin{eqnarray}\label{alg}
(e(x) \otimes g)( e(y) \otimes h) &=& e(x)e(g.y) \otimes gh \notag \\
&=& \delta_{x, g.y}e(x) \otimes gh
\end{eqnarray}
for $x, y \in X$ and $g, h \in G$.
One readily checks that this is \emph{associative}. There is a decomposition into
$2$-sided ideals
\begin{eqnarray}\label{dirsum}
D_X(G) = \bigoplus_Y D_Y(G)
\end{eqnarray}
where $Y$ ranges over the (transitive) $G$-orbits of $X$. For the most part, this reduces questions about $D_X(G)$ for general $X$ to the transitive case.

\bigskip
A special example of this construction is the quantum double of $G$. Here, we take $X = G_{conj}$, i.e. $X=G$ and the left action of $G$ is left conjugation $g: x \mapsto gxg^{-1}$. In this case we write $D(G)$ in place of  $D_{G_{conj}}(G)$.  $D(G)$ is in fact a Hopf algebra, but at the moment we only require the algebra structure.

\bigskip
Next we describe the  category of (left-) $D_X(G)$-modules
(cf. \cite{DPR}, \cite{M1}, \cite{KMM}). For $x \in X$ and a left Stab$_G(x)$-module $V$, set $V_x = e(x) \otimes V$. This is a left
$\mathbb{C}[X]^* \otimes$ Stab$_G(x)$-module via
\begin{eqnarray}\label{raction}
(e(y) \otimes g). (e(x) \otimes v) = e(y)e(x) \otimes g.v = \delta_{x, y} e(x) \otimes g.v, \ g \in \mbox{Stab}_G(x).
\end{eqnarray}
From (\ref{alg}) it follows that $\mathbb{C}[X]^* \otimes$ Stab$_G(x)$ is a left ideal in $D_X(G)$, 
so that \linebreak
$D_X(G) \otimes_{\mathbb{C}[X]^* \otimes Stab_G(x)} V_x$ is a left $D_X(G)$-module.

\begin{propn}\label{propn2.1} Suppose that $X$ is a \emph{transitive} $G$-set and $x \in X$. Then the map 
\begin{eqnarray}\label{Morequiv}
&& \mathbb{C}[\mbox{Stab}_G(x)]\mbox{-Mod} \rightarrow D_X(G)\mbox{-Mod}, \notag \\
&& \hspace{2.7cm} V \mapsto D_X(G) \otimes_{\mathbb{C}[X]^* \otimes Stab_G(x)} V_x,
\end{eqnarray}
is a Morita equivalence.
\end{propn}
\begin{pf} For this and more, see \cite{M1} and Section 3 of \cite{KMM}. In these references 
$X$ is a group, but this is not necessary and the proofs to go through without change. $\hfill \Box$
\end{pf}

\begin{propn}\label{propn2.2} Suppose that $X$ is a transitive $G$-set and $x \in X$. The following hold:
\begin{eqnarray*}
&&(a) \ \mbox{If $p$ does not divide $|X|$ then $(D_X(G),$ Stab$_G(x))$ is an M-pair.}\\
&&(b) \ \mbox{If $p$ divides $|X|$ then $\mu(D_X(G)) = 0$.}
\end{eqnarray*}
\end{propn}
\begin{pf} By Proposition \ref{propn2.1} and (\ref{Morequiv}), the simple left $D_X(G)$-modules
are precisely the modules $D_X(G) \otimes_{\mathbb{C}[G]^* \otimes Stab_G(x)} V_x$ as $V$ ranges over the simple left modules
for Stab$_G(x)$.  

\bigskip
Let $T$ be a set of right coset representatives in $G$ for Stab$_G(x)$, so that
there is a disjoint union $G = \cup_{t \in T} t\mbox{Stab}_G(x)$. If $t \in T,s \in$ Stab$_G(x)$ and
$y \in X$ then $e(y) \otimes ts = (e(y) \otimes t)(e(t^{-1}.y) \otimes s$. Using (\ref{raction}) we have
\begin{eqnarray*}
D_X(G) \otimes_{\mathbb{C}[G]^* \otimes Stab_G(x)} V_x &=& \sum_{s, t, y} (e(y) \otimes t)(e(t^{-1}.y)
\otimes s) e(x) \otimes V \\
&=&\sum_{t, y} \delta_{x, t^{-1}.y} (e(y) \otimes t)e(x) \otimes V \\
&=&\sum_{t}  (e(t.x) \otimes t)e(x) \otimes V,
\end{eqnarray*}
the last sum being a \emph{direct sum.} Because $X$ is transitive, it follows that
\begin{eqnarray*}
\dim (D_X(G) \otimes_{\mathbb{C}[G]^* \otimes Stab_G(x)} V_x) 
&=& |T|\dim V \\
&=& |X| \dim V.
 \end{eqnarray*}
 
 It follows that in the Morita equivalence (\ref{Morequiv}), modules of dimension $d$ are mapped to modules of dimension $d|X|$. Parts (a) and (b) follow immediately from this. $\hfill \Box$
\end{pf}

\begin{lemma}\label{lemma2.3}
$\mu(D_X(G)) = \sum_Y \mu(D_Y(G))$ where $Y$ ranges over the $G$-orbits of $X$.
\end{lemma}
\begin{pf} This follows from the decomposition (\ref{dirsum}) into $2$-sided ideals. $\hfill \Box$
\end{pf}

\bigskip
We will also need the following standard result. 
\begin{lemma}\label{lemma2.4}  The number of $G$-orbits of $X$ of cardinality coprime to $p$ is equal to the number of $N$-orbits of $X$ of cardinality coprime to $p$.\end{lemma}
\begin{pf} By considering the decomposition of $X$ into $G$-orbits, we see that we must prove the following assertion:
\begin{eqnarray}\label{transassert}
&&\mbox{If $X$ is a \emph{transitive} $G$-set, then either (a) $p$ divides $|X|$ and $N$ has} \notag\\
&&\mbox{\emph{no} orbits of cardinality coprime to $p$, or (b) $p$ does \emph{not} divide $|X|$}\\
&&\mbox{and there is a \emph{unique} $N$-orbit of cardinality coprime to $p$}.\notag
\end{eqnarray}

Assume, then, that $X$ is a transitive $G$-set, and that $x, y \in X$ lie in $N$-orbits of cardinality coprime to $p$. In such an $N$-orbit, $P$ must fix at least one, and therefore all,
elements in the $N$-orbit. In particular, $P$ lies in the stabilizers of both $x$ and $y$. 
By transitivity there is $g \in G$ with $g.x = y$. Then $P$ and $P^g$ are both Sylow $p$-subgroups
of Stab$_G(x)$ and by Sylow's theorem  there is $t \in$ Stab$_G(x)$ such that $P^{gt} = P$. Then
$gt \in N$ and $(gt).x = y$. This shows that $x, y$ lie in the \emph{same} $N$-orbit, so that there is at most one $N$-orbit of cardinality coprime to $p$. (\ref{transassert}) is easily deduced from this, and the Lemma is proved. $\hfill \Box$
\end{pf}

\bigskip
Consider the following statements:
\begin{eqnarray*}
&&\mbox{MC$\mathbf{1}$}: (\mathbb{C}[G],\mathbb{C}[N]) \ \mbox{is an M-pair for all $G$},\\
&&\mbox{MCD}: (D(G), D(N)) \ \mbox{is an M-pair for all $G$},\\
&&\mbox{MCX}: (D_X(G), D_X(N)) \ \mbox{is an M-pair for all $G$ and all $X$},\\
&&\mbox{MCT}:  (D_X(G), D_X(N)) \ \mbox{is an M-pair for all $G$ and all \emph{transitive} $X$}.
\end{eqnarray*}
The McKay Conjecture is of course the assertion that MC$\mathbf{1}$ is true.

\begin{thm}\label{thm2.5} $MC\mathbf{1}, MCD, MCX$ and $MCT$ are \emph{equivalent} statements.
\end{thm}
\begin{pf} $MCX \Leftrightarrow MCT$: This follows from Lemma \ref{lemma2.4}. 

\medskip
\noindent
$MCT \Rightarrow MC\mathbf{1}$:  This holds because if $X = \mathbf{1}$ is the one-element set,
then $D_X(G) = \mathbb{C}[G]$. 

\medskip
\noindent
$MC\mathbf{1} \Rightarrow MCT$: 
Let $X$ be a transitive $G$-set. If $p$ divides $|X|$ then $\mu(D_X(G)) = 0$ by Proposition \ref{propn2.2}(b). Similarly,
$\mu(D_X(N))=0$ by Lemma \ref{lemma2.3}, (\ref{transassert})(a), and Proposition \ref{propn2.2}(b) (applied to $N$).  Now suppose that $p$ does \emph{not} divide $|X|$. Then 
$\mu(D_X(G)) = \mu(\mathbb{C}[\mbox{Stab}_G(x)])$ for any $x \in X$ by Proposition \ref{propn2.2}(a).
Moreover by Lemma \ref{lemma2.4} there is a \emph{unique} $N$-orbit of cardinality coprime to $p$,
call it
 $Y \subseteq X$. By Lemmas \ref{lemma2.3}, (\ref{transassert}) and Proposition \ref{propn2.2} once more we find that $\mu(D_X(N)) =\mu(D_Y(N)) = \mu(\mathbb{C}[\mbox{Stab}_N(y)])$ for $y \in Y$.
Note that because $|Y|$ is coprime to $p$ then $P \leq$ Stab$_G(y)$ and
$N_{Stab_G(y)}(P) =$ Stab$_N(y)$. The assumption that $MC\mathbf{1}$ holds (applied to
Stab$_G(y)$) tells us that $\mu(\mathbb{C}[\mbox{Stab}_G(y)])
= \mu(\mathbb{C}[\mbox{Stab}_N(y)])$, whence $\mu(D_X(G)) = \mu(D_X(N))$.

\medskip
\noindent
$MCX \Rightarrow MCD$:  Let $Y_1, \hdots , Y_h$ be the $N$-orbits
of $G_{conj}$ of cardinality coprime to $p$. We have $Y_i \subseteq N$ for each index $i$, so that they are also the $N$-orbits of $N_{conj}$ of cardinality coprime to $p$. Taking $X = G_{conj}$, MCX together with 
Lemma \ref{lemma2.3} and Proposition \ref{propn2.2}, we conclude that
$\mu(D(G)) = \mu(D_{G_{conj}}(N)) = \sum_{i=1}^h \mu(D_{Y_i}(N)) = \mu(D_{N_{conj}}(N))
= \mu(D(N)),$ as required. 

\medskip
\noindent
 $MCD \Rightarrow MC\mathbf{1}$: We prove this using induction on $|G|$. Retain the notation of the last paragraph, and choose $y_i \in Y_i$. By Lemma \ref{lemma2.4},
$y_1, \hdots, y_h$ are representatives for the $G$-orbits of $G_{conj}$ (ie. conjugacy classes
of $G$) of cardinality coprime to $p$. By Proposition \ref{propn2.2}, $\mu(D(G)) = \sum_{i=1}^h \mu(\mathbb{C}[C_G(y_i)])$.
Since each $y_i \in N$, we similarly have $\mu(D(N)) = \sum_{i=1}^h \mu(\mathbb{C}[C_N(y_i)])$. If 
$C_G(y_i)$ is a \emph{proper} subgroup of $G$ then by induction 
$\mu(\mathbb{C}[C_G(y_i)]) = \mu(\mathbb{C}[C_N(y_i)])$. Then the assumption $MCD$ tells us that
$\sum_{i'} \mu(\mathbb{C}[G]) = 
\sum_{i'} \mu(\mathbb{C}[N])$ where $i'$ ranges over those indices for which
$y_{i'}$ lies in the \emph{center} $Z(G)$ of $G$. We conclude that
$|Z(G)| \mu(\mathbb{C}[G])= |Z(G)|\mu(\mathbb{C}[N])$, whence $\mu(\mathbb{C}[G])= \mu(\mathbb{C}[N])$. This completes the proof of the Theorem. $\hfill \Box$
 \end{pf}

\section{Twisted Algebras}
In this Section we explain how to extend the results of the previous Section to the 
\emph{twisted case}, i.e. the incorporation of a cocycle. Let $\theta \in Z^2(G, \mathbb{C}^*)$ be a (normalized) multiplicative $2$-cocycle. Thus $\theta: G^2 \rightarrow \mathbb{C}^*$ satisfies the identities
\begin{eqnarray*}
\theta(h, k)\theta(g, hk) &=& \theta(gh, k)\theta(g, h), \ \ g, h, k \in G, \\
\theta(1, g) &=& \theta (g, 1) = 1.
\end{eqnarray*}
The corresponding
twisted group algebra is $\mathbb{C}^{\theta}[G]$. It has the same underlying linear space as
$\mathbb{C}[G]$ with multiplication
$g \circ h = \theta(g, h)gh$ for $g, h \in G$. The cocycle identities ensure that this is an associative algebra with identity element $1$. For a subgroup $H \leq G$ we identify $\theta$ with its
\emph{restriction} Res$^G_H \theta$ to $H$. Then $\mathbb{C}^{\theta}[H]$ is a subalgebra
of $\mathbb{C}^{\theta}[G]$. For more information on this subject, including results that we use below, see for example \cite{CR}.

\bigskip
The cohomological analog of Proposition \ref{propn2.2}(b) is the following
\begin{lemma}\label{lemma3.1} Suppose that the \emph{cohomology class} $[\theta] \in H^2(G, \mathbb{C}^*)$
determined by $\theta$ has order $k$. If $k$ is divisible by $p$ then
$\mu(\mathbb{C}^{\theta}[G]) = \mu(\mathbb{C}^{\theta}[N]) = 0$.
\end{lemma}
\begin{pf} One knows (loc. cit.) that  there is a central extension
\begin{eqnarray*}
1 \rightarrow Z  \rightarrow L \stackrel{\pi}{\rightarrow} G \rightarrow 1
\end{eqnarray*}
such that $\mathbb{Z}_k \cong Z \leq L',$ and $\mathbb{C}^{\theta}[G]$ is the algebra summand of $\mathbb{C}[L]$ corresponding to the 
irreducible representations of $L$ in which a generator $z$ of $Z$ acts as multiplication by a prescribed primitive $k$th. root of unity, say $\lambda$. If $V$ is a simple $\mathbb{C}^{\theta}[G]$-module of dimension $d$ then the determinant of $z$ considered as operator on $V$ is clearly $\lambda^d$.
On the other hand $z \in L'$, so that this determinant is necessarily $1$. So $\lambda^d = 1$,
whence $k|d$. In particular, if $p|k$ then $\mu(\mathbb{C}^{\theta}[G]) = 0$.

\bigskip
Now it is well-known that the restriction map $\Res^G_N: H^2(G, \mathbb{C}^*) \rightarrow
H^2(N, \mathbb{C}^*)$ is an \emph{injection} on the $p$-part of $H^2(G, \mathbb{C}^*)$.
In particular, if $p|k$ then $\Res^G_N [\theta]$ has divisible by $p$. Then the result of the last paragraph also applies to $\mathbb{C}^{\theta}[N]$, and we obtain $\mu(\mathbb{C}^{\theta}[N]) = 0$. This completes the proof of the Lemma. $\hfill \Box$
\end{pf}

\bigskip
The McKay Conjecture implies that the twisted analog is also true. This is the content of
\begin{propn} Suppose that MC$\mathbf{1}$ holds. Then  $(\mathbb{C}^{\theta}[G]), \mathbb{C}^{\theta}[N])$
is an M-pair
for all $G$ and all $\theta$.
\end{propn}
\begin{pf} Let the notation be as in Lemma \ref{lemma3.1}. Although it is not really necessary to do so, because of Lemma \ref{lemma3.1} we may, and shall, assume that $k$ is not divisible by $p$. 
Let $P_1$ be a Sylow $p$-subgroup of $L$ with $\pi: P_1 \stackrel{\cong}{\rightarrow} P$. Applying MC$\mathbf{1}$ to  pairs
$(L/Z_0, N_L(P_1)/Z_0)$ with  $Z_0 \leq Z$, we see that
the number $l$ of irreducible representations of both $\mathbb{C}[L]$ and $\mathbb{C}[N_L(P_1)]$ which have degree coprime to $p$ and in which $z$ acts as \emph{some} primitive $k$th. root of unity are equal. Since $\mathbb{C}^{\theta}[G]$
is the algebra summand of $\mathbb{C}[L]$ corresponding to $\lambda$, then 
$\mu(\mathbb{C}^{\theta}[G]) = l/\phi(k)$. On the other hand, $N_L(P_1) = ZK$ where
$K \leq L$ is such that $\mathbb{C}^{\theta}[N]$ is the algebra summand   of
$\mathbb{C}[K]$ corresponding to $\lambda^t$, where $t|k$ and $k/t$ is the order of
$\Res^G_N[\theta]$.   By slightly modifying the previous argument, we also find that 
$\mu(\mathbb{C}^{\theta}[N]) = l/\phi(k),$  and the Proposition is proved. $\hfill \Box$
\end{pf}

\bigskip
We can now treat the twisted version of $D_X(G)$.
Let $U = U(\mathbb{C}[X]^*)$ be the group of units in $\mathbb{C}[X]^*$. Then
\begin{eqnarray*}
U = \left\{ \sum \lambda_x e(x) \ | \ \lambda_x \not= 0 \right\}
\end{eqnarray*}
is a \emph{multiplicative} left $G$-module. Let $\alpha \in Z^2(G, U)$ be a normalized $2$-cocycle with coefficients in $U$, and set $\alpha(g, h) = \sum_{x \in X}
\alpha_x(g, h)e(x)$. Here, the cocycle property amounts to the identity
\begin{eqnarray}\label{alphamult}
\alpha_x(g, h)\alpha_x(gh, k) = \alpha_x(g, hk)\alpha_{g^{-1}.x}(h, k).
\end{eqnarray}

\medskip
Define $D_X^{\alpha}(G)$ to be the linear space $D_X(G)$ with multiplication being the twisted 
version of (\ref{alg}). That is, 
\begin{eqnarray}\label{tp}
 (e(x)\otimes g)(e(y) \otimes h) = \alpha_x(g, h) \delta_{x, g.y} e(x) \otimes gh.
 \end{eqnarray}
(\ref{alphamult}) is exactly what is needed to show that (\ref{tp}) is \emph{associative}.
Note also from (\ref{alphamult}) that for fixed $x \in X, \ \alpha_x$ defines an element
in $Z^2(\mbox{Stab}_G(x), \mathbb{C}^*)$ and that  as a subspace of $D_X^{\alpha}(G), \ S(x) \cong \mathbb{C}^{\alpha_x}[\mbox{Stab}_G(x)].$  The proof of Proposition \ref{propn2.1} still applies in this situation (cf. \cite{KMM}). It provides a Morita equivalence of categories
\begin{eqnarray*}
\mathbb{C}^{\alpha_x}[\mbox{Stab}_G(x)]\mbox{-Mod} \stackrel{\sim}{\rightarrow}D^{\alpha}_X(G)
\mbox{-Mod}.
\end{eqnarray*}
The proof of the twisted version of Theorem \ref{thm2.5} then goes through too. We just state a part of this as  
\begin{thm}\label{thm3.3} Suppose that MC$\mathbf{1}$ holds. Then $(D_X^{\alpha}(G), (D_X^{\alpha}(N))$ is an M-pair for all $G, X$ and $\alpha$. $\hfill \Box$
\end{thm}

\bigskip
Special cases of $D^{\alpha}_X(G)$ include certain kinds of crossed products and abelian extensions
of Hopf algebras. See, for example, \cite{KMM} for further details.

\bigskip
Once again  the case of the quantum double, when $X = G_{conj}$,  is of special interest (cf. \cite{CK}, \cite{DM}, \cite{DPR}, \cite{M1} for more details and further background.) Here, one twists $D(G)$ by a
normalized \emph{three cocycle} $\omega \in Z^3(G, \mathbb{C}^*)$. The resulting object is denoted by
$D^{\omega}(G)$. It is a quasi-Hopf algebra, but not a Hopf algebra in general. To connect with previous paragraphs, we observe that there is
a map (\cite{DPR})
\begin{eqnarray*}
Z^3(G, \mathbb{C}^*) \rightarrow Z^2(G, G_{conj}) 
\end{eqnarray*}
for which
\begin{eqnarray}\label{thetas}
\alpha_x(g, h) = \frac{\omega(x, g, h)\omega(g, h, (gh)^{-1}x(gh))}{\omega(g, g^{-1}xg, h)}.
\end{eqnarray}
 There is a natural interpretation of this map in terms of  the loop space $LBG$, but we will not need it. The twisted product in $D^{\omega}(G)$
is as in (\ref{tp}) using (\ref{thetas}). This gives the algebra structure, and as before leads to
\begin{thm}\label{thm3.4} Suppose that MC$\mathbf{1}$ holds. 
Then $(D^{\omega}(G), D^{\omega}(N))$ is an M-pair for all $G$ and $\omega$. $\hfill \Box$
\end{thm}

\bigskip
The statement and proof of Theorem \ref{thm3.4} only requires the algebra structure of
$D^{\omega}(G)$. However, we will make use of other structural features of $D^{\omega}(G)$
in the next Section.

\section{Orbifolds and Ribbon Categories}
We refer the interested reader to \cite{DM} for background concerning veretx operator algebras. Let $V$ be a holomorphic vertex operator algebra admitting $G$ as a group of automorphisms, with $V^G$  the subalgebra of $G$-invariants. One expects that the module category $V^G$-Mod is a (braided, ribbon) tensor category and that it is equivalent to the tensor category $D^{\omega}(G)$-Mod for a $3$-cocycle $\omega$ which describes the associativity constraint in $V^G$-Mod.  If this is so, we deduce from Theorem \ref{thm3.4} that there are bijections between the simple objects of $V^G$-Mod and $V^N$-Mod which themselves correspond to the 
simple modules of $D^{\omega}(G)$-Mod and $D^{\omega}(N)$-Mod respectively which have dimension coprime to $p$. 

\medskip
We seek a direct definition of an M-pair for modules over
orbifolds such as $V^G$ and $V^N$. We cannot use (\ref{Mpairdef}) as it stands because it makes no sense for infinite-dimensional spaces such as a module over a vertex operator algebra. Instead, we can make use of the expected structure of $V^G$-Mod as a ribbon tensor category, whereby the objects have a \emph{quantum dimension}. Indeed, $D^{\omega}(G)$-Mod has a \emph{canonical} ribbon structure (cf. \cite{AC}, \cite{MN}), and  the quantum dimensions of simple objects  are the usual dimensions. Granted the equivalence of $V^G$-Mod and $D^{\omega}(G)$-Mod, it follows that
the quantum dimension of simple objects in $V^G$-Mod are also integers. Then the definition of an M-pair makes sense if we use quantum dimension in place of dimension. 

\medskip
Thus we arrive at the following situation:  a pair of ribbon fusion categories $\mathcal{G}, \mathcal{N}$ whose simple objects have quantum dimensions that are rational integers. We  say that 
$(\mathcal{G}, \mathcal{N})$ is an M-pair if $\mu(\mathcal{G}) = \mu(\mathcal{N})$, where we use (\ref{Mpairdef}) with quantum dimension in place of dimension in order to define $\mu$. As we have seen, taking $\mathcal{G}$ to be $\mathbb{C}[G]$-Mod, $D^{\omega}(G)$-Mod or $V^G$-Mod and $\mathcal{N}$ to be $\mathbb{C}[N]$-Mod, $D^{\omega}(N)$-Mod or $V^N$-Mod respectively
results (conjecturally) in an M-pair. Furthermore, the three versions of MC for groups, quantum doubles of groups, and holomorphic orbifolds, are \emph{equivalent}.

\end{document}